\documentclass[12pt,journal,draftcls,letterpaper,onecolumn]{IEEEtran}

\usepackage{lscape}
\usepackage{cite}
\usepackage{amsmath}
\usepackage{amsmath}
\usepackage{amsfonts}   
\usepackage{amssymb}    
\usepackage{float}
\usepackage{subfigure}
\usepackage{float}
\usepackage{amsmath}
\usepackage{color}
\usepackage{graphicx}
\usepackage{epsfig}
\usepackage{algorithm}
\usepackage[noend]{algorithmic}
\usepackage[doublespacing]{setspace}
\begin{document}
\title{Resource Allocation via Sum-Rate Maximization in the Uplink of Multi-Cell OFDMA Networks}
\author{\IEEEauthorblockN {Hina Tabassum$^*$, Zaher Dawy$^{**}$, and Mohamed Slim Alouini$^*$}\\
\IEEEauthorblockA{$^*$Physical Sciences and Engineering Division, King Abdullah University of Science and Technology (KAUST), Thuwal, Mekkah Province, Saudi Arabia \\
\IEEEauthorblockA{$^{**}$Department of Electrical and Computer Engineering, American University of Beirut (AUB), Beirut, Lebanon} \\
Email: hina.tabassum@kaust.edu.sa, zd03@aub.edu.lb, slim.alouini@kaust.edu.sa}}
\maketitle
\begin{abstract}
In this paper, we consider  maximizing the sum-rate in the uplink of a multi-cell OFDMA network. The problem has a non-convex combinatorial structure and is known to be NP hard. Due to the inherent complexity of implementing the optimal solution, firstly, we derive an upper and lower bound to the optimal average network throughput. Moreover, we investigate the performance of a near optimal single cell resource allocation scheme in the presence of ICI which leads to another easily computable lower bound. We then develop a centralized sub-optimal  scheme that is composed of a geometric programming based power control phase in conjunction with an iterative subcarrier allocation phase. Although, the scheme is computationally complex, it provides an effective benchmark for low complexity schemes even without the power control phase. 
Finally, we propose less complex centralized and distributed schemes that are well-suited for practical scenarios. The computational complexity of all schemes is analyzed and performance is compared through simulations. Simulation results demonstrate that the proposed low complexity schemes can achieve comparable performance to the centralized sub-optimal scheme in various scenarios. Moreover, comparisons with the upper and lower bounds provide insight on the performance gap between the proposed schemes and the optimal solution.

\end{abstract}
\IEEEpeerreviewmaketitle
\section{Introduction}
Dynamic resource allocation plays a central role in the air interface design of state-of-the-art OFDMA-based cellular technologies. In this paper, we focus our attention on maximizing the overall network throughput by optimizing the allocation of resources (i.e., subcarriers and powers) jointly in a multi-cell uplink OFDMA network. The goal is to develop efficient resource allocation schemes that takes into account the  inter-cell interference (ICI) while considering universal frequency reuse. The solution of such problem is difficult to achieve optimally due to its NP hard combinatorial nature and  high dimensionality. 

The sum rate maximization problem is extensively studied for the downlink in OFDMA networks. The optimal strategy in the downlink is to separately optimize subcarrier and power allocation, i.e., allocate a subcarrier to the user with best channel and then perform water-filling over the allocated subcarriers \cite{Downlink}. However, the problem becomes more challenging in the uplink scenario due to the individual power constraint at each user. Simply allocating a subcarrier to the user with best channel quality may affect the network performance considerably, as some active users may have better channel gains but low transmission powers on a specific subcarrier. 


Most of the recent work in the context of multi-cell OFDMA networks \cite{6,7}, aims at minimizing the overall transmitted power, i.e., linear objective with pre-defined rate constraints. In \cite{8}, the authors investigated scaling laws for upper and lower bounds of the downlink capacity in the asymptotic regime. Furthermore, in some recent literature, low complexity distributed game theoretic solutions are also studied. However, the schemes are iterative and optimality is not guaranteed \cite{9}. An auction based approach is discussed in \cite{10}, where the authors proposed a joint auction and dual decomposition based technique for the resource allocation problem. The technique is asymptotically optimal as the number of subcarriers in every cell goes to infinity. However, this may not be true for finite number of carriers. In summary, all these approaches are sub-optimal and no criteria are used to calibrate their performance gap with respect to the optimal solution. 

Motivated by the above discussed facts, we consider the problem of optimized resource allocation in the uplink of multi-cell OFDMA networks.
Firstly, we compute an upper bound (UB) and lower bound (LB) to the optimal average network throughput. Also, we study the effect of ICI on the performance of the near-optimal single cell resource allocation scheme proposed in \cite{kim} which leads to another simple lower bound. Simulation results show that this lower bound is slightly loose but can be computed easily. Since the computation of the optimal solution is exhaustive, we then propose a centralized sub-optimal resource allocation scheme which uses a geometric programming (GP) based power control phase in conjunction with a heuristic subcarrier allocation phase. The proposed scheme possesses an iterative and computationally intensive subcarrier allocation phase. However, it can serve as an effective benchmark for the less complex schemes even without the power control phase. Furthermore, the power control phase is discussed in this paper for both high and general signal to interference plus noise ratio (SINR) regimes. Finally, we propose and evaluate less complex centralized and distributed schemes that are suitable for practical implementation. 

The rest of the paper is organized as follows: In Section II, the system model is defined and problem is formulated.  In Section III, the bounds are derived and their complexity is analyzed. In Section IV and Section V, the proposed centralized and distributed  schemes are explained. Section VI presents numerical results followed by concluding remarks in Section VII.\\ 
\textbf{Notation}: Throughout the paper, we denote the sets of real and complex vectors of $N$ elements by  $\mathbb{R}^N$ and $\mathbb{C}^N$, respectively. Matrices are represented using boldface upper case letters while bold face lower case letters are used for vectors. $\mathcal{N}(0,\sigma^2)$ represents a zero mean Gaussian random variable with variance $\sigma^2$.

\section{System Model and Problem Formulation}
A network of $L$ cells with a set of $K$ users in each cell $l$ is considered.  Full reuse of the spectrum is assumed in all the cells (i.e., frequency reuse =1). Each base station (BS) is assumed to have $N$ orthogonal subcarriers, and each subcarrier can be allocated to a single user per cell. The average network throughput $\mathcal{C}$ is a function of both subcarrier and power allocation variables. The sum rate maximization problem is formulated as follows using the standard Shannon capacity formula, 
$\mathcal{C}_{n,k,l}=\mathrm{log_2}(1+\gamma_{n,k,l})$, where $\mathcal{C}_{n,k,l}$ and $\gamma_{n,k,l}$ represent the throughput and SINR of the $k^{\mathrm{th}}$  user at $n^{\mathrm{th}}$  subcarrier in cell $l$, respectively: 
\begin{equation}
\label{obj}
		 \underset{p_{n,k,l},\alpha_{n,k,l}}{\text{maximize}}\:\:\:\:
			\sum\limits_{l=1}^L 
				\sum\limits_{k=1}^{K} 
					\sum\limits_{n=1}^N 
					\alpha_{n,k,l}
					\log_2
					\left(1+
						\frac{\textstyle p_{n,k,l} h_{n,k,l}}
								{\textstyle \sigma^2+I_{n,l}}
						\right)\\
\end{equation}
\begin{equation}
\label{powconst}
	\begin{aligned}
\text{subject to}	& & 
			\sum_{n=1}^N p_{n,k,l}\leq \mathrm{P}_{k,\mathrm{max}},\;\:\forall k,\forall l\:\:\:\:\:\:\:\:\:\:\:\:\:\:\:\:\:\:\:\:\:\:\:\:\:\:\:\:
\end{aligned}
\end{equation}
\begin{equation}
\label{allconst}
	\begin{aligned}
			\sum_{k=1}^{K}\alpha_{n,k,l}=1,\;\:\forall n,\forall l \:\:\:\:\:\:\:\:\:\:\:\:\:\:\:\:\:\:\:\:\\
\end{aligned}
\end{equation}
\begin{equation}
\label{allintconst}
	\begin{aligned}
		\alpha_{n,k,l}\in\{0,1\},\;\:\forall n,\forall l,\forall k\:\:\:\:\:\:\:\:\:\:\:\:\:\:\:\:\:\:\:\:\:\:\\
\end{aligned}
\end{equation}
In \eqref{obj}, $ I_{n,l}= \sum_{j=1,j\neq l}^L \sum_{k=1}^{K} \alpha_{n,k,j}  p_{n,k,j} g_{n,k,jl}$ represents the cumulative interference at $n^{\mathrm{th}}$ subcarrier in cell $l$ from the users in all other cells, $p_{n,k,l}$ denotes the power transmitted by $k^{\mathrm{th}}$ user at the  $ n^{\mathrm{th}}$  subcarrier in cell $l$, $\alpha_{n,k,l}$ represents the allocation of $k^{\mathrm{th}}$  user at the $n^{\mathrm{th}}$ subcarrier in cell $l$ and $h_{n,k,l}$ is the channel gain of $k^{\mathrm{th}}$  user at the $n^{\mathrm{th}}$  subcarrier in cell $l$. Constraint \eqref{powconst} implies that the power spent by $k^{\mathrm{th}}$ user on its allocated subcarriers cannot exceed the maximum available power, $P_{k_{\mathrm{max}}}$. For each cell, we collect the power allocation variables $p_{n,k,l}$ in a vector  ${\bf{p}}_{n,l}=[p_{n,1,l},p_{n,2,l},....,p_{n,K,l}] $ and then stack all the vectors in a power matrix ${\bf{P}}_l$ of cell $l$ where $ {\bf{P}}_{l}\:\in\: \mathbb{R}^{N\times K}$. Constraint \eqref{allconst} restricts the allocation of each subcarrier to only one user. The channel gains $h_{n,k,l}$ and binary allocation variables $\alpha_{n,k,l}$ are stacked up similarly in the matrices $ {\bf{H}}_{l}$ and $ {\bf{A}}_{l}$, respectively, where ${\bf{A}}_l, {\bf{H}}_l\:\in\: \mathbb{R}^{N\times K} $. Moreover, we define $g_{n,k,lj}$ as the interfering gain from the $k^{\mathrm{th}}$ user in cell $l$ to cell $j$, $\forall j\neq l$ at $n^{\mathrm{th}}$ subcarrier. We collect these interfering gains into a vector  ${\bf{g}}_{n,lj}=[g_{n,1,lj},g_{n,2,lj}  ....,g_{n,K,lj}] $ and then stack all the vectors in a matrix ${\bf{G}}_{lj} \:\in\: \mathbb{R}^{N\times K}$.

\subsection{Optimal Problem Formulation in High SINR Regime}
Assuming perfect knowledge of channel gains at a centralized controller, the optimal solution for \eqref{obj} can be computed in the high SINR regime by an exhaustive search over all possible combinations of the allocations. For each possible allocation, optimum powers can be computed by transforming \eqref{obj} into a GP. Note that the power allocation problem is in itself a known non-convex problem for the general SINR regime \cite{11}. However, in the high SINR regime the problem becomes a convex GP problem. For a given set of allocation variables and considering a high SINR regime, the objective function in \eqref{obj} can be rewritten as follows:
\begin{equation}
\label{HighSNR1}
\begin{aligned}
		& \underset{{p_{n,k,l}}}{\text{maximize}}	& & 
			\sum_{l=1}^L \sum_{k=1}^{K} \sum_{n=1}^N 
			\alpha_{n,k,l} \mathrm{log_2} \left (\frac {p_{n,k,l}h_{n,k,l}}{\sigma^2+ I_{n,l}}\right)
\end{aligned}
\end{equation}
Maximizing the SINRs is equivalent to minimizing the interference to signal ratio:
\begin{equation}
\label{HighSNR2}
\begin{aligned}
		& \underset{p_{n,k,l}}{\text{minimize}}	& & 
 \sum_{l=1}^L \sum_{k=1}^{K} \sum_{n=1}^N \alpha_{n,k,l}\mathrm{log_2} \left(\frac {\sigma^2+ I_{n,l}}{p_{n,k,l} h_{n,k,l}}\right)
\end{aligned}
\end{equation} 
Equivalently, \eqref{obj} can be reformulated for high SINR regime and given allocation variables as follows:
\begin{equation}
\label{HighSNR3}
	\begin{aligned}
		& \underset{p_{n,k,l}}{\text{minimize}}	& & 
			\log_2
			\prod\limits_{l=1}^L 
				\prod\limits_{k=1}^{K} 
					\prod\limits_{n=1}^N 
					\left (\frac {\sigma^2+ I_{n,l}}{p_{n,k,l} h_{n,k,l}}\right)^{\alpha_{n,k,l}}
						\\
		& \text{subject to}	& &
			\sum_{n=1}^N \alpha_{n,k,l} p_{n,k,l}\leq\mathrm{P}_{k,\mathrm{max}},\;\:\forall k,\forall l\\
\end{aligned}
\end{equation}
Note that the numerator in \eqref{HighSNR3} is a posynomial and the denominator is a monomial, hence \eqref{HighSNR3} is a GP problem in standard form that can be solved optimally through efficient interior point methods \cite{13} after performing the logarithmic transformation of variables \cite{11}. However, even for small dimensions, it is not recommendable to compute the optimal solution, due to the  huge computational complexity $O(K^{LN})$ associated with the exhaustive search based subcarrier allocation phase. In addition, the GP based power allocation method discussed above has two restrictions: high-SINR assumption and centralized time-consuming computations. Due to the mentioned facts, there is a need to develop bounds and sub-optimal resource allocation schemes for multi-cell OFDMA networks.
\section{Bounds on the Network Throughput}
\subsection{Lower Bound on the Optimal Network Throughput} 
A LB for the optimum multi-cell network throughput can be computed by considering worst case ICI. Observing the dependency of ICI on the subcarrier allocation and power allocation variables, we assume that each user in each cell is transmitting on each subcarrier with its maximum power. A simple LB for the average network throughput $\mathcal{C}$ taking the worst case ICI into account can be written as follows:
\begin{equation}
\label{LB}
\mathcal{C}({\bf{A}}_l,{\bf{P}}_l)\geq \frac{1}{L}  \sum_{l=1}^L \sum_{k=1}^{K} \sum_{n=1}^N \alpha_{n,k,l} \mathrm{log_2} \left(1+\frac {p_{n,k,l} h_{n,k,l}}{\sigma^2+ \xi_{n,l}}\right)
\end{equation}
where $\xi_{n,l}= \sum_{j=1,j\neq l}^L \sum_{k=1}^{K} P_{k,\mathrm{max}}g_{n,k,jl}$. 

A tighter LB can be derived by
using Algorithm 1 where each subcarrier is allocated to the user that maximizes ${Q}_{n,k,l}$ where:
\begin{equation}
Q_{n,k,l}=\frac{p_{n,k,l}h_{n,k,l}}{\xi_{n,l}+\sigma^2}
\end{equation}  
Thus, $Q_{n,k,l}$ is an SINR term for each user $k$ at each subcarrier $n$ in each cell $l$ assuming worst case interference. We collect these SINR terms into a vector  ${\bf{q}}_{n,l}=[q_{n,1,l},q_{n,2,l}  ....,q_{n,K,l}] $ and then stack all the vectors in a matrix ${\bf{Q}}_{l} \:\in\: \mathbb{R}^{N\times K}$. The resulting allocations based on this criteria are then used to compute the LB network throughput using \eqref{obj}. 

Note that if $\xi_{n,l}=0$, than ${Q}_{n,k,l}$ becomes the marginal rate which is shown to be a near-optimal criterion in single cell network scenarios without ICI \cite{kim}. Moreover, equal power allocation has insignificant performance loss in high SINR regime compared to the optimal water-filling solution \cite{kim,yaoma}, thus power equalization is implemented in Algorithm 1. For the low SINR regime, we can incorporate water-filling  rather than equalization in a straightforward manner.
\begin{algorithm}
\caption{Computing LB and UB Allocations in Cell $l$}         
\label{alg1}   
\begin{enumerate}
\item \textbf{Input:}$[{\bf{H}}_l],[{\bf{A}}_l],[{\bf{P}}_l],[{\bf{G}}_{jl}]$ where $\alpha_{n,k,l}=0, p_{n,k,l}=P_{k,\mathrm{max}}/N \:\: \forall k,\forall n $
\item For each user $k$ in cell $l$, power is divided equally over all of its allocated subcarriers and the remaining unallocated subcarriers of the system.
\item Using $[{\bf{P}}_l ]$ from step 2, $[{\bf{H}}_l ]$ and $[{\bf{G}}_{jl}]$, compute the matrix ${\bf{Q}}_l$ for each cell $l$. 
\item Find the $(n,k)$ pair that has the maximum value of ${Q}_{n,k,l}$. Allocate subcarrier $n$ to user $k$. 
\item Delete the $n^{\mathrm{th}}$ subcarrier from the set of unallocated subcarriers.\\ {\textbf{If}} there are still unallocated subcarriers in the system go to step 2.\\ {\textbf{else}} terminate after distributing the maximum power equally at each user over all of its assigned subcarriers.
\end{enumerate}
\vskip -3pt
\end{algorithm}

 \subsection{Upper Bound on the Optimal Network Throughput}
Establishing an UB is significantly important in order to calibrate the performance of sub-optimal resource allocation schemes with respect to the optimal solution. The UB can be derived by ignoring the effect of ICI in all the cells. This can be achieved by substituting $ \xi_{n,l}=0$ in Algorithm 1, i.e., ${Q}_{n,k,l}=\frac{p_{n,k,l}h_{n,k,l}}{\sigma^2}$: 
\begin{equation}
\label{UB}
\mathcal{C}({\bf{A}}_l,{\bf{P}}_l) \leq \frac{1}{L}  \sum_{l=1}^L \sum_{k=1}^{K} \sum_{n=1}^N \alpha_{n,k,l} \mathrm {log_2} \left (1+\frac {p_{n,k,l} h_{n,k,l}}{\sigma^2}\right)
\end{equation}
The allocations computed by Algorithm 1 are near optimal since they are based on a criterion which is shown to be near optimal in the context of single cell scenarios\cite{kim,ng,yaoma}. The average network throughput revealed by these allocations could be highly optimistic for multi-cell scenarios. Thus, we can investigate the impact of ICI by simply computing the throughput using \eqref{obj} instead of \eqref{UB} with these allocations. Computing throughput in this way helps to analyze the degradation in the performance when the single cell near-optimal allocations are used in multi-cell network scenarios with ICI.
\subsection{Complexity Analysis}
The $(n,k)$ pair at which the term $Q_{n,k,l}$ becomes maximum is allocated (Step 4), which has a complexity of a two dimensional search, i.e., $O(KN)$. However, as soon as a subcarrier is assigned, each user updates its power as defined in Algorithm 1. This process iterates until all the subcarriers in all the cells are allocated and, thus, the time complexity of Algorithm 1 is $O(KN^{2})$.
\subsection{A Motivating Example}
Consider an example with two cells, two users and two subcarriers. Each user can transmit with a maximum power of 1 W. Assume ${\bf{H}}_1$= [1 0.9; 0.8 0.7] and ${\bf{H}}_2$= [1 0.9; 0.8 0.7]. Single cell allocation strategies that aim to maximize the local throughput of each cell suggest ${\bf{A}}_1,{\bf{P}}_1$ and ${\bf{A}}_2,{\bf{P}}_2 =$[1 0; 0 1]. Computing the UB using \eqref{UB} results in 1.7655 bps/Hz/cell where $\sigma^2=1$. Now, assume the knowledge of interfering link gains at each BS, i.e., ${\bf{G}}_{12}$= [0.9 0.2; 0.2 0.9] and ${\bf{G}}_{21}$= [0.7 0.1; 0.1 0.7]. Computing the throughput again while keeping the single cell allocations and taking into account the interfering gains leads to an average network throughput of 1.1137 bps/Hz/cell. However, better allocations are possible if we consider ${\bf{A}}_1,{\bf{P}}_1$ and ${\bf{A}}_2,{\bf{P}}_2 =$ [0 1; 1 0] as per the criterion discussed in Section IV which enhances the resulting average network throughput to 1.5977 bps/Hz/cell.
 \section{Sub-Optimal Centralized Resource Allocation Schemes}
Considering the high intricacy of implementing the optimal solution, we develop a two-stage centralized scheme. 
In comparison to the centralized scheme presented in \cite{hina}, the subcarrier allocation phase of the developed scheme is iterative and computationally intensive. However, the performance is better even without the power allocation phase and, thus, it can provide an effective benchmark for low complexity schemes. Compared to \cite{hina}, we also discuss the power allocation phase for the general SINR regime.

\subsection{Centralized Scheme A}
In the proposed scheme, we split the resource allocation procedure into two phases: subcarrier allocation phase and power allocation phase. It is important to note that the subcarrier allocation phase involves a power equalization step, thus, it is not totally independent of power allocation.\\
\textbf{Phase I: Subcarrier Allocation}
\begin{itemize}
\item \textit{Initial Allocation}: Firstly, we define the term for the allocation of resources to the users as follows:
\begin{equation}
\label{ki}
\chi_{n,k,l} =\frac{p_{n,k,l}h_{n,k,l}}{\sum_{j=1,j \neq l}^L P_{k,\mathrm{max}}{g}_{n,k,lj}}
\end{equation} 
This criterion guarantees the selection of the users who possess not only better power-gain product but also they offer less interference to the neighbor cells. The denominator ${\sum_{j=1,j \neq l}^L P_{k,\mathrm{max}}{g}_{n,k,lj}}$ accounts for the maximum aggregate interference that the $k^{\mathrm{th}}$ user in cell $l$, may cause to all cells. Even though this criterion is heuristic, it improves the performance compared to the traditional C/I scheme (which gives nearly similar results as our lower bound). Once the initial allocations are computed, we can calculate the initial throughput of the network $\mathcal{C}_o$ using \eqref{obj}.
\item \textit{Maximize Throughput Iteratively until Convergence}: 
In this step, we select any cell $l$ and subcarrier $n$ arbitrarily and re-perform the allocation at this subcarrier considering the other cell allocations fixed, i.e., $I_{n,l}$ remains fixed. More explicitly, we compute $\mathcal{C}_{n,k,l}=\mathrm{log_2}(1+\frac{\textstyle p_{n,k,l} h_{n,k,l}}{\textstyle \sigma^2+I_{n,l}})$ for all users in cell $l$ one by one and select the user which gives the maximum incremental throughput at subcarrier $n$, i.e., $\mathcal{C}_{n,k,l}-\mathcal{C}_o$. 
Note that, in order to compute $\mathcal{C}_{n,k,l}$, we need to compute $p_{n,k,l}$ which can be obtained simply by dividing $P_{k,\mathrm{max}}$ equally among all the fixed allocated subcarriers of user $k$ and the new one which is currently under observation. 

Once the reallocation is done at subcarrier $n$, we move to the next subcarrier in cell $l$  and so on. As the new allocations are computed for cell $l$, we calculate the new increased network throughput $\mathcal{C}_{\mathrm{new}}$ and move to another cell $j$. The whole process is repeated again with $\mathcal{C}_o=\mathcal{C}_{\mathrm{new}}$ until convergence to a desired accuracy is achieved.
\end{itemize}
\textbf{Phase II: Power Allocation}\\
Once the subcarrier allocation is done, the optimal powers can then be calculated for the high SINR regime or for the general SINR regime through solving a series of GPs using successive convex approximation which is a provably convergent heuristic \cite{11}. This approach is known to compute globally optimal power allocations in many cases. Thus, for given allocations, \eqref{obj} can be formulated for the general SINR regime as follows:
\begin{equation}
\label{GENSINR}
	\begin{aligned}
		& \underset{p_{n,k,l}}{\text{minimize}}	& & 
			\log_2
			\prod\limits_{l=1}^L 
				\prod\limits_{k=1}^{K} 
					\prod\limits_{n=1}^N 
					\left (\frac {\sigma^2+ I_{n,l}}{p_{n,k,l} h_{n,k,l}+\sigma^2+ I_{n,l}}\right)^{\alpha_{n,k,l}}
						\\
		& \text{subject to}	& &
			\sum_{n=1}^N \alpha_{n,k,l}p_{n,k,l}\leq\mathrm{P}_{k,\mathrm{max}},\;\:\forall k,\forall l\\
\end{aligned}
\end{equation}
Note that the numerator and denominator in \eqref{GENSINR} are posynomials and minimizing a ratio between two posynomials is referred to be a truly non-convex NP hard intractable problem known as complimentary GP. However, this problem can be transformed into GP by letting the denominator $f(p)={p_{n,k,l} h_{n,k,l}+\sigma^2+ I_{n,l}}=\sum_{l=1}^{L}\sum_{k=1}^{K}u_{n,k,l}(p)$ and approximating the denominator $f(p)$ with a monomial using the arithmetic/geometric mean inequality as follows:
\begin{equation}
\label{AMGM}
\sum_{l=1}^{L}\sum_{k=1}^{K}u_{n,k,l}(p) \geq \prod_{l=1}^{L}\prod_{k=1}^{K} {\left(\frac{u_{n,k,l}(p)}{s_{n,k,l}}\right)}^{s_{n,k,l}}
\end{equation}
where $s_{n,k,l}=\frac{u_{n,k,l}(p_0)}{f(p_0)}$. Thus, the problem can be solved by extending the single condensation method presented in \cite{11} for multi-cell scenario. The details of centralized scheme A are presented in Algorithm 2.
\begin{algorithm}
\caption{\: Centralized Scheme A }        
\label{alg2}
\begin{enumerate}
\item \textbf{Input:}$[{\bf{H}}_l],[{\bf{A}}_l],[{\bf{P}}_l],[{\bf{G}}_{lj}]$ where $\alpha_{n,k,l}=0, p_{n,k,l}=P_{k,\mathrm{max}}/N \:\: \forall k,\forall n $\\
{\textbf{Subcarrier Allocation (Phase I)}} \\ 
\textit{Initial Allocation:} 
\item For each user $k$ in cell $l$, power is divided equally over all of its allocated subcarriers and the remaining unallocated subcarriers of the system.
\item Using $[{\bf{P}}_l ]$ from step 2, $[{\bf{H}}_l ]$ and $[{\bf{G}}_{lj}]$,  compute $\chi_{n,k,l}$ for every $k^{\mathrm{th}}$ user at $n^{\mathrm{th}}$ subcarrier in cell $l$.
\item Find the $(n,k)$ pair that has the maximum value of ${\chi}_{n,k,l}$. Allocate subcarrier $n$ to user $k$. 
\item Delete the $n^{\mathrm{th}}$ subcarrier from the set of unallocated subcarriers.\\ {\textbf{If}} there are still unallocated subcarriers in the system go to step 2,\\ {\textbf{else}} terminate after distributing the maximum power at each user over all of its assigned subcarriers 
\item Compute $\mathcal{C}_o$ \\
\textit{Maximize Throughput Iteratively until Convergence}\\
{\textbf {do while}}($\mathcal{C}_{\mathrm{new}}-\mathcal{C}_o \geq \epsilon)$\\
$l=1$, {\textbf {do while}} $l \leq L, l=l+1 $\\
\:\:\: $n=1$, {\textbf {do while}} $n \leq N, n=n+1 $\\
\:\:\:\:\:$k=1$, {\textbf {do while}} $k \leq K, k=k+1 $
\item Allocate the subcarrier $n$ to user $k$.
\item Compute $p_{n,k,l}$ by dividing $P_{k,\mathrm{max}}$ equally among the allocated subcarriers.
\item Compute $\mathcal{C}_{n,k,l}-\mathcal{C}_o$\\
\:\:\:\:\:{\textbf{end}}
\item Allocate subcarrier $n$ to the user who maximizes $\mathcal{C}_{n,k,l}-\mathcal{C}_o$\\
\:\:\:{\textbf{end}}
\item Compute $\mathcal{C}_{\mathrm{new}}$ using \eqref{obj}. 
\item $C_o=C_{\mathrm{new}}$\\
{\textbf{end}}\\
{\textbf{Power Allocation (Phase II)}} 
\item Compute the optimal powers ${\bf{P}}_l$ in the high SINR regime \eqref{HighSNR3} given the allocations from Phase I.
\item For general SINR regime, take ${\bf{P}}_l$ from step 13 as an initial starting point.
\item Using  ${\bf{P}}_l$, evaluate  $p_{n,k,l} h_{n,k,l}+\sigma^2+ I_{n,l}$ for each allocated user $k$ in cell $l$ at subcarrier $n$.
\item Compute the weights $s_{n,k,l}$ as follows:\\
$s_{n,k,l}=\frac{u_{n,k,l}}{f(p)}$
\item	Approximate the posynomial using \eqref{AMGM}.
\item Solve the approximated GP using any available commercial software \cite{13}
\item Go to step 15 using ${\bf{P}}_l$ of step 18 until convergence.
\end{enumerate}
\end{algorithm}
\subsection{Centralized Scheme A: Complexity Analysis}
The initial allocation phase has a complexity of $O(KN^2)$ which is the same as 
Algorithm 1. Next, we perform a one dimensional search for the user in cell $l$ with maximum incremental throughput at subcarrier $n$. The process is repeated for each subcarrier and cell. Thus, the computational complexity of this step is $O(KNL)$. Since, the process continues until convergence, (i.e., $M$ iterations), the complexity of this step can be written as $O(KNLM)$. Finally, the total complexity of subcarrier allocation phase is $O(KN^2 +NKLM)$. 

The complexity of Phase II is difficult to determine, however, it can be measured in terms of the degree of difficulty (DoD) that in turn relies on the number of constraints and variables associated with the GP \cite{12}. Since we are dealing with $LK$ power constraints and $LKN$ power variables, the total computational complexity of centralized scheme A is $O(KN^2 +NKLM)+DoD(LKN)$. Apparently it seems that implementing centralized GP/successive GP based schemes may not be a good choice for practical implementations. However, in order to reduce the complexity and DoD of the power control phase, we have developed the following less complex centralized scheme.
\subsection{Centralized Scheme B}
In this scheme, firstly the subcarriers are allocated in each cell $l$ using the heuristic criterion defined in \eqref{ki}. The allocation of each subcarrier is followed by the power allocation phase (based on equalization) as mentioned in the initial allocation phase of Algorithm 2 (i.e., Steps 1 to 5). Once the subcarrier allocations are finalized,
we then compute GP based powers  for the allocated users at any arbitrarily selected subcarrier $n$ in all cells. Setting the equalization based powers $p_{{n,k,l}_{\mathrm{eq}}}$ as the upper bound on $p_{n,k,l}$ and considering a high SINR regime, we now define the following less complex GP problem with the objective to maximize the throughput at the $n^{\mathrm{th}}$ subcarrier: 
\begin{equation}
\label{centB}
	\begin{aligned}
		& \underset{p_{n,k,l}}{\text{minimize}}	& & 
			\log_2
			\prod\limits_{l=1}^L 
					\left (\frac {\sigma^2+ I_{n,l}}{p_{n,k,l} h_{n,k,l}}\right)^{\alpha_{n,k,l}}
						\\
		& \text{subject to}	& &
			 p_{n,k,l}\leq {p}_{{n,k,l}_{\mathrm{eq}}},\;\: \forall l
\end{aligned}
\end{equation}
Clearly, the resulting GP based power of each competing user at subcarrier $n$ in the different cells may not succeed in achieving the upper bound, due to the ICI effect. We call this power as left-over power. The left-over power can then be distributed equally among the remaining allocated subcarriers of the user.
The procedure is detailed in Algorithm~3.

Since at the end of the initial allocation phase, the subcarrier allocations become fixed and the total power is distributed equally among the allocated subcarriers of a user, we cannot set an upper bound which depicts higher power than the previously allocated power. If we do so, this may cause power reduction or even no power at some other allocated subcarrier of that user in order to maintain the total power constraint. Thus, this may results in an invalid subcarrier allocation.

\begin{algorithm}
\caption{\: Centralized Scheme B }        
\label{alg3}   
\begin{enumerate}
\item Repeat Steps 1 to 5 of Algorithm 2, i.e., initial allocation phase.\\
\:\:\: $n=1$, {\textbf {do while}} $n \leq N, n=n+1 $
\item Compute the GP based powers ${{p}}_{n,k,l}$ of the allocated users at any subcarrier $n$ considering a  high SINR regime using \eqref{centB}.
\item For each user allocated in a cell $l$ at any subcarrier $n$, divide the left-over power equally among the remaining allocated subcarriers of the user.
\item Remove the subcarrier $n$ from the set of unallocated subcarriers.
{\textbf{end}}
\end{enumerate}
\end{algorithm}

Next follows an example which  demonstrates the significance of GP as well as centralized scheme B over equal power allocation.  Consider ${\bf{H}}_1, {\bf{H}}_2$= $[0.30 \:\:\:\:\:  0.25;\:\: 0.04 \:\:\:\:\: 0.15]\times 10^{-9}$, ${\bf{G}}_{12} = [0.06 \:\:\:\:\: 0.05; \:\:0.16  \:\:\:\:\:   0.06]\times 10^{-11}$ and $ {\bf{G}}_{21}$=$[0.14 \:\:\:\:\:  0.69;\:\: 0.76 \:\:\:\:\: 0.1935]\times 10^{-11}$. The equal power allocations dictate ${\bf{P}}_1=[0\:\:\:\:\:   0.5; \:\:0\:\:\:\:\:  0.5]$ and ${\bf{P}}_2=[0.5 \:\:\:\:\:  0;\:\: 0.5 \:\:\:\:\:  0]$ which leads to an average network throughput of  11.8392 bps/Hz/cell. However, computing the GP based powers results in ${\bf{P}}_1=[0 \:\:\:\:\:   0.53;  \:\:
   0   \:\:\:\:\:  0.47]$ and ${\bf{P}}_2=[0.38 \:\:\:\:\:    0;
   \:\:  0.62 \:\:\:\:\: 0]$ which lead to a maximum average network throughput of 17.2734 bps/Hz/cell.

\subsection{Centralized Scheme B: Complexity Analysis} 
The initial allocation phase has a complexity of $O(KN^2)$ which is the same as Algorithm~1. Since \eqref{centB} has $L$ constraints and variables, the complexity of the power control phase is significantly reduced. Although this procedure restricts the degree of freedom offered by GP, numerical results show that the network throughput remains comparable with reduced complexity.

\section{Distributed Resource Allocation Scheme}
In the centralized strategy, we assume that $\chi_{n,k,l}$ is known, i.e., every BS knows the interfering gains offered by its users to the neighboring BSs. The interfering gains are based on path loss, shadowing and fading. Assuming the knowledge of local user positions at each BS, the path loss of local users toward the first tier of interfering cells can be determined, however, the knowledge of shadowing and fading gains is difficult to assume in practical scenarios. Thus, in the distributed approach, we compute our results without using the knowledge of shadowing and fading interfering gains.

Each BS performs the subcarrier allocations without taking ICI into account. In other words we compute single cell near optimal allocations using Algorithm 1. The allocation decisions are locally made at each BS and do not need collaboration. Once the allocations are decided, each cell shares them with all other interfering cells. The GP based optimal powers in \eqref{HighSNR3} can then be evaluated in a distributed way using dual decomposition methods by first performing the log transformation of the variables, i.e., $\mathrm{ln} p_{n,k,l}= 
\tilde{p}_{n,k,l}$ and  $\mathrm{ln} p_{n,k,j}= \tilde{p}_{n,k,j}$, then adding auxiliary variable $\mathrm{ln} z_{n,lj}=\tilde{z}_{n,lj}$ where ${z_{n,lj}}={p_{n,k,j}}$ in order to transfer the coupling in the objective to coupling in the constraints \cite{11}. For given allocations, the problem in \eqref{HighSNR3} can thus be written in a distributed way as follows:

\begin{equation}
\label{dist1}
	\begin{aligned}
		\underset{\tilde{z}_{n,lj},\tilde{p}_{n,k,l}}{\text{minimize}}	& 
				\sum\limits_{l=1}^{L} 
					\sum\limits_{n=1}^N 
						\log_2
					\left (\frac {\sigma^2+  \sum_{j=1,j \neq l}^L  g_{n,k,jl} e^{\tilde{z}_{n,lj}}}{e^{\tilde{p}_{n,k,l}}h_{n,k,l}}\right)\\	
		\text{subject to} &  \sum_{n=1}^N  e^{\tilde{p}_{n,k,l}}\leq\mathrm{P}_{k,\mathrm{max}},\;\:\forall k,\forall l
\\
     & {z}_{n,lj} = {p}_{n,k,j},\;\:\forall n,\forall l
\end{aligned}
\end{equation}
Since the computational complexity of \eqref{dist1} is high as it has $LK$ power constraints and $LKN$ variables,
we present the dual decomposition of \eqref{centB}  which is more suitable for practical scenarios and has a lower computational complexity. Moreover,  the objective function in \eqref{dist1} not only depends on the powers of local users $p_{n,k,l}$ but also on the power of users sharing the same subcarrier in neighboring cells $p_{n,k,j}$. Thus, in order to minimize the objective in \eqref{dist1}, each BS requires the knowledge of interfering gains and interfering transmit powers, that may lead to significant overhead to exchange control information. Thus, in order to obtain a practical distributed solution, we keep a local copy of each of the effective received powers i.e., $z_{n,lj}=g_{n,k,jl}{p}_{n,k,j}$ \cite{11}. \eqref{centB} can then be formulated in a distributed way as follows:
\begin{equation}
\label{dist2} 
	\begin{aligned}
		\underset{\tilde{z}_{n,lj},\tilde{p}_{n,k,l}}{\text{minimize}}	& \sum\limits_{l=1}^L 
						\log_2
					\left (\frac {\sigma^2+ \sum_{j=1,j \neq l}^L e^{\tilde{z}_{n,lj}}}{e^{\tilde{p}_{n,k,l}}h_{n,k,l}}\right)\\	
		\text{subject to} & \:\:\:\:e^{\tilde{p}_{n,k,l}}\leq\mathrm{P}_{n,k,l,\mathrm{eq}},\;\:\forall l
\\
     & {\tilde{z}_{n,lj}} =  \tilde{g}_{n,k,jl} +{\tilde{p}_{n,k,j}}
\end{aligned}
\end{equation} 
The partial lagrange $\mathcal{L}(\tilde{p}_{n,k,l},\tilde{z}_{n,lj},\lambda_{l},\eta_{n,lj})$ for  \eqref{dist2} can then be written explicitly as follows:
\begin{equation}
\label{Lagrange}
\sum\limits_{l=1}^L 
						\log_2
					\left (\frac {\sigma^2+ \sum_{j=1,j \neq l}^L e^{\tilde{z}_{n,lj}}}{e^{\tilde{p}_{n,k,l}}h_{n,k,l}}\right)+ \sum\limits_{l=1}^L \sum\limits_{j=1,j \neq l}^L \eta_{n,lj}\left({\tilde{z}_{n,lj}} -\tilde{g}_{n,k,jl}-{\tilde{p}_{n,k,j}}\right)+\sum_{l=1}^L \lambda_l\left(
e^{\tilde{p}_{n,k,l}}-\mathrm{P}_{n,k,l,\mathrm{eq}}\right)
\end{equation}
Eq. \eqref{Lagrange} can be decomposed into $L$ sub-problems with local variables $\tilde{p}_{n,k,l},\tilde{z}_{n,lj},\lambda_{l}$ and coupling variable $\eta_{n,lj}$. The simple lagrangian $\mathcal{L}_l$ for each cell $l$ can then be written as follows:
\begin{equation}
\label{locaLag}
					\mathcal{L}_l=	\log_2
					\left(\frac {\sigma^2+ \sum_{j=1,j \neq l}^L e^{\tilde{z}_{n,lj}}}{e^{\tilde{p}_{n,k,l}}h_{n,k,l}}\right)+ 
\sum\limits_{j=1,j \neq l}^L \eta_{n,lj}{\tilde{z}_{n,lj}} -\left(\sum\limits_{j=1,j \neq l}^L \eta_{n,jl}\right){\tilde{p}_{n,k,l}} + \lambda_l\left(
{\tilde{p}_{n,k,l}}-\mathrm{P}_{n,k,l,\mathrm{eq}}\right)
\end{equation}
where $\lambda_{l}$ is the lagrange multiplier for the inequality constraints and $\eta_{n,lj}$ are the consistency prices. Thus, the minimization of \eqref{locaLag} with respect to the local variables can be done in a distributed way at all BSs. At every iteration, each cell $l$ receives the term $\left(\sum\limits_{j=1,j \neq l}^L \eta_{n,jl}\right)$ by message
passing and minimizes the local Lagrangian \eqref{locaLag} with
respect to  the local variables $\tilde{p}_{n,k,l},\tilde{z}_{n,lj}, \lambda_{l}$ subject to the local constraints. In order to obtain $\eta_{n,lj}$ the following master lagrange dual problem has to be solved:
\begin{equation}
\label{dist}
	\begin{aligned}
		\underset{{\eta_{n,lj}}}{\text{maximize}} \sum_{l=1}^L \underset{\tilde{p}_{n,k,l},\tilde{z}_{n,lj},\lambda_l} {\text{minimize} }\:\:\:\mathcal{L}_l
\end{aligned}
\end{equation}
A simple way to solve \eqref{dist} is to use the following subgradient update for
the consistency prices:
\begin{equation}
\label{update}
\eta_{n,lj}(t+1)=\eta_{n,lj}(t)+ ({\delta}/t)\left(\tilde{z}_{n,lj}-\mathrm{log_2}\: {p}_{n,k,j} {g}_{n,k,jl} \right)
\end{equation}
In summary, each BS minimizes \eqref{locaLag} in parallel with respect to the local variables after receiving the term $\sum_{j=1,j\neq l}^L  \eta_{n,jl}$. Each BS then estimates the received interference ${z}_{n,lj}$ from each cell and update the local consistency prices using \eqref{update}. Finally, each BS broadcast them by message passing to all BSs.
Note that $\delta$ in \eqref{update} represents the step size and is non-negative.

%
\section{Performance Evaluation}
A cellular OFDMA network is considered where the radius of each cell is 
assumed to be 1~km. The maximum user transmit power is considered to be 1 W. The channel gain is defined as follows:
\begin{eqnarray*} 
h_{n,k,l}=(-122-10 \gamma\: {\mathrm {log_{10}}} d_{k,l})-\mathcal{N}(0,\sigma^2 )+10\: {\mathrm{log}}_{10}F_{n,k,l}  \:\:\:\:\: d> d_{\mathrm{ref}} \\
h_{n,k,l}=(-122-10 \gamma\: {\mathrm {log_{10}}} d_{\mathrm{ref}})-\mathcal{N}(0,\sigma^2 )+10\: {\mathrm{log}}_{10}F_{n,k,l}  \:\:\:\:\: d<d_{\mathrm{ref}}
\end{eqnarray*}
where $d_{\mathrm{ref}}$ is the reference distance and is set equal to 0.05 km, $d_{k,l}$ is the distance of the $k^{\mathrm{th}}$ user from the $l^{\mathrm{th}}$ BS. The first term denotes the path loss where $\gamma$ is the path loss exponent and is set equal to 3. The second term represents log-normal shadowing with a mean of 0 dB and a standard deviation of 8 dB. The last factor $F_{n,k,l} $ corresponds to Rayleigh fading. The bandwidth of the system is
assumed to be 20 MHz with a noise power spectral density of
$8.6455 \times 10^{-15}$~W/Hz at each receiver. The channel conditions are assumed to be fixed during a frame. The interfering gains from the $j^{\mathrm{th}}$ interfering cell to the cell of interest $l$ are computed as follows:
$$g_{n,k,jl}=(-122-10\:\gamma \:\mathrm{log_{10}} d_{k,l})-\mathcal{N}(0,\sigma^2 )+10\:\mathrm{log}_{10}F_{n,k,j}$$
where $d_{k,l}$ is the distance between the $k^{\mathrm{th}}$ user in the interfering cell $j$ and the $l^{\mathrm{th}}$ BS. We consider the following two simulation scenarios:
\begin{itemize}
\item \textbf{Scenario A:} Users are equidistant from the BS and placed at equally spaced angles.
\item \textbf{Scenario B:} Users are assumed to be uniformly distributed across the whole cellular area.
\end{itemize} 

In Table 1, we compare the performance and complexity of the centralized and distributed schemes with the derived bounds and the optimal solution in high SINR regime. The optimal solution is computed by the exhaustive search based subcarrier allocation phase detailed in Section II. All users are placed at equal distance $d$ from the BS and at equally spaced angles (i.e., scenario A).  The results are taken after averaging over 100 channel realizations. The simulation results show that the performance gap between the benchmark centralized scheme A (with power control) and the optimal solution is negligible compared to the low complexity centralized and distributed schemes.
However, this observation may not remain valid for bigger network scenarios. Moreover, as $d$ increases the degradation of the average network throughput is evident. 

In Fig. 1, Fig. 2 and Fig. 3, we present the performance of the centralized scheme A, centralized scheme B  and the distributed scheme for two cells, four cells and seven cells, respectively. The results have been taken after averaging over 10,000 channel realizations and are shown for both simulation scenarios. The performance of all schemes is calibrated using the established upper and lower bounds.
Since the centralized scheme A has computationally exhaustive power allocation phase, the results are presented for the subcarrier allocation phase of Algorithm 2 only. However, it can be observed that the scheme still has the capability to serve as a suitable benchmark for the developed low complexity schemes.
In order to highlight the significance of the less complex GP problem defined in \eqref{centB}, we also present the performance of the centralized scheme B without power control.  

For the two cell scenario, the performance gap between the centralized schemes is negligible and they give nearly similar results. However, as the number of cells increases the performance gain of the centralized scheme A is evident over all schemes even without power control. Moreover, it is also important to note the significant degradation in the performance of centralized scheme B without power control phase. This degradation is found to be increasing with the increase in number of cells. It is also worth to mention here that the proposed less complex GP problem \eqref{centB} can be implemented in a distributed way using the techniques explained in \cite{11} and can be used with any set of subcarrier allocations. Thus, in the distributed approach we use the near optimal single cell allocations in conjunction with the less complex GP problem \eqref{centB}. The significance of the power control phase can be observed easily from the results which becomes more evident for high number of cells .

Moreover, the presented results depict the reduction in the average network throughput as the number of interfering cells increases. The performance gap of the proposed schemes increases with respect to the evaluated UB. Even though the UB is not tight and reflects an over optimistic average network throughput, it provides an idea on the performance gap between the proposed schemes and the optimal solution. 

In Fig. 4, we assume that the users in each cell are placed at equally spaced angles from 0 to $2\pi$. The performance evaluation of all schemes has been done by changing the user positions from cell center to the edge of the cell. It is observed that the performance gap increases between the centralized and distributed schemes as users approach the boundaries of the cell. 

\section{Conclusion}
In this paper, we developed an upper bound and a lower bound to the optimal average network throughput in multi-cell uplink OFDMA networks. We also investigated the severe effects of ICI on the performance of a single cell near optimal resource allocation scheme. Moreover, we proposed a benchmark centralized scheme which is useful to study the performance gap of the low complexity centralized and distributed resource allocation schemes developed later with respect to the optimal solution. All schemes are compared to the exhaustive search based optimal solution and derived upper and lower bounds for various scenarios. The schemes are evaluated and compared in terms of network throughput and computational complexity.

\section{Acknowledgement}
The authors would like to express their sincere gratitude to Dr. Elias Yaccoub for useful discussions toward the success of this work.
\bibliography{IEEEfull,references}
\bibliographystyle{IEEEtran}

\newpage
\begin{table}[H]
\renewcommand{\arraystretch}{1.9}
\begin{center}
\caption{Average network throughput (in bps/Hz/cell) of the derived bounds, centralized and distributed schemes for $L$=2 cells and $N=6$ subcarriers/cell}
\resizebox{1\textwidth}{!}
{\begin{tabular}{|l||l|l||l|l||l|l||l|}
\hline\hline
&\multicolumn{2}{l|}{$K$=2 Users}&\multicolumn{2}{l|}{$K$=4 Users}&\multicolumn{2}{l|}{$K$=6 Users}& Computation Complexity\\
\cline{2-7}
&$d$=0.5km& $d$=0.9km&$d$=0.5km&$d$=0.9km&$d$=0.5km&$d$=0.9km\\
\hline\hline
UB          &44.2642      &33.1294      &55.7414   & 42.8390   &60.6901 &49.6214& $O(KN^2)$ \\
Optimal     &37.1168       &29.8642      &47.9975     & 35.6520   &52.1299  &41.0121 & $O(K^{NL}) +DoD(LKN)$\\
Centralized A &36.8061     &28.6973     &46.4765    & 34.0713   &51.2868 &40.5845& $O(KN^2 +NKLM)+DoD(LKN)$\\
Centralized B &36.4755     &27.0352      &45.6239     & 33.4280   &49.7971   &38.7237 & $O(KN^2)+DoD(L)$\\
Distributed &35.3623   &25.9976    &43.5918    & 31.9231   &48.8887 & 38.0050&$O(KN^2)+DoD(L)$\\
LB          &35.0966      &25.8635     &42.5509      &31.0261   & 48.1571
& 37.7996 &$O(KN^2)$\\
\hline 
\end{tabular}}
\end{center}
\end{table} 
\newpage
\begin{figure}[H]
  \centering
  \includegraphics[width=5in]{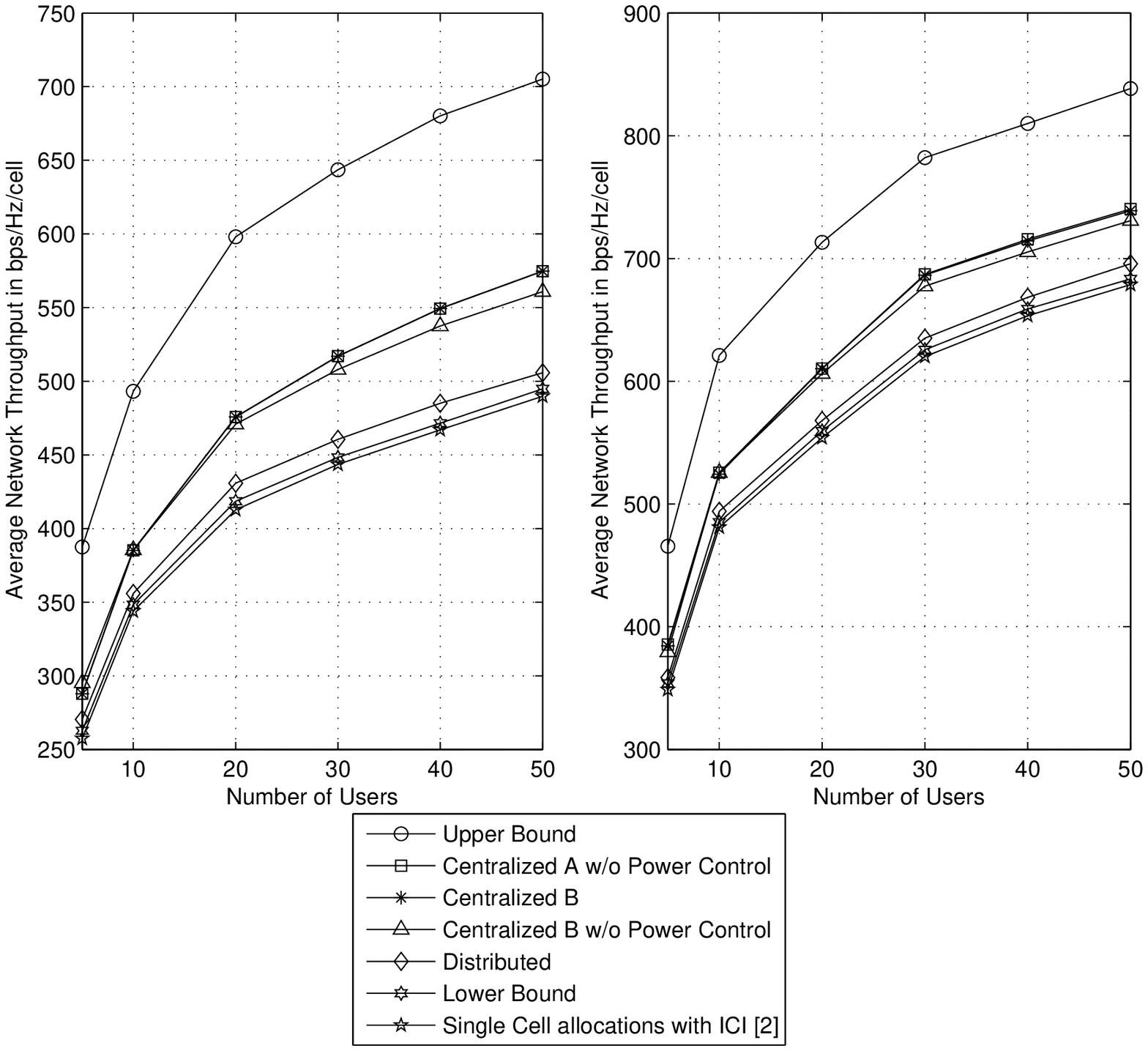}
  \caption{Comparison of all proposed schemes for $L$=2 cells, (a) Scenario A (b) Scenario B: Users are placed at 0.9 km from~ BS}
\label{Cap}
\end{figure}
\newpage
\begin{figure}[H]
  \centering
  \includegraphics[width=5in]{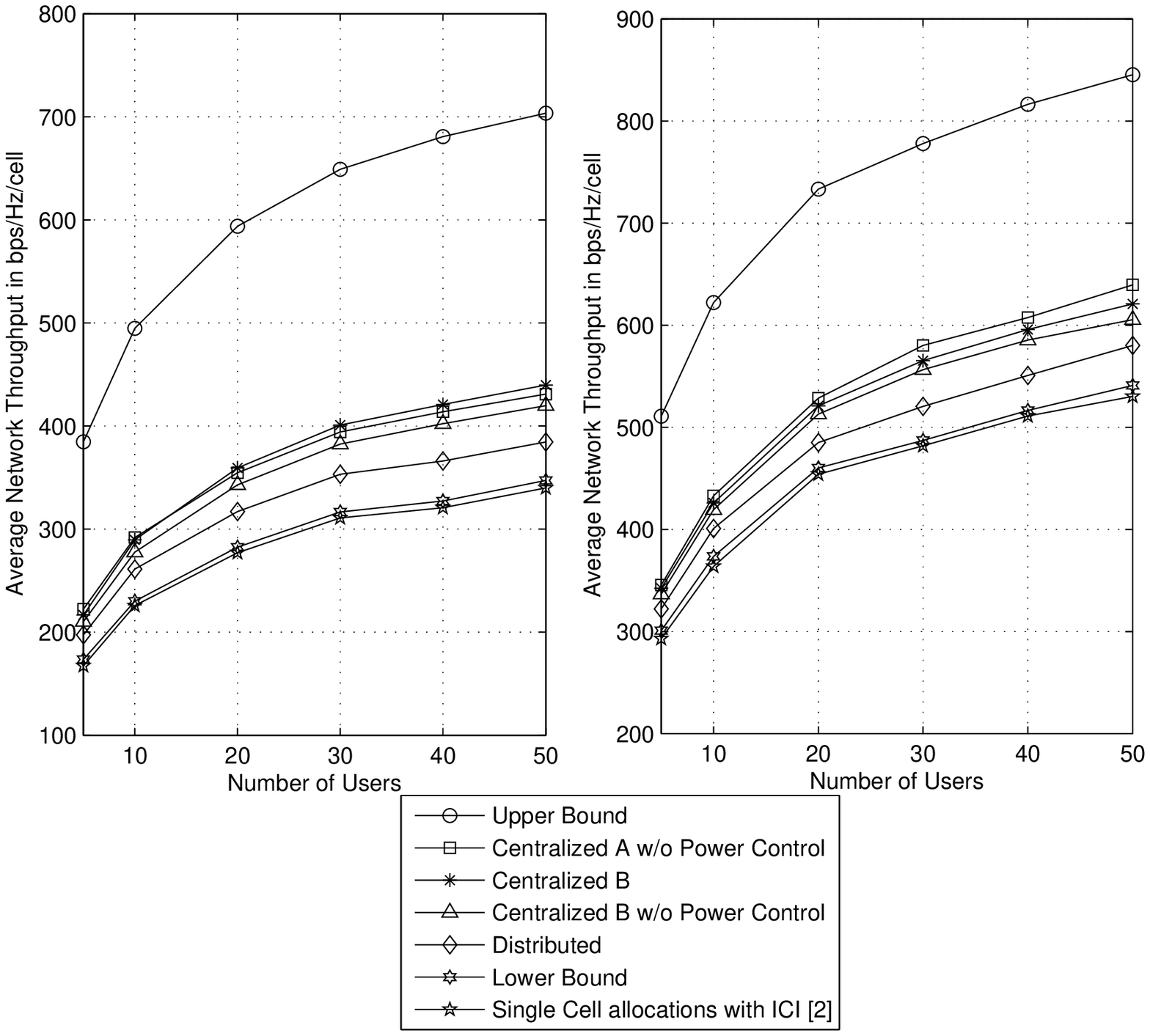}
  \caption{Comparison of all proposed schemes for $L$=4 cells, (a) Scenario A (b) Scenario B: Users are placed at 0.9 km from~ BS}
\label{Cap}
\end{figure}
\newpage
\begin{figure}[H]
  \centering
  \includegraphics[width=5in]{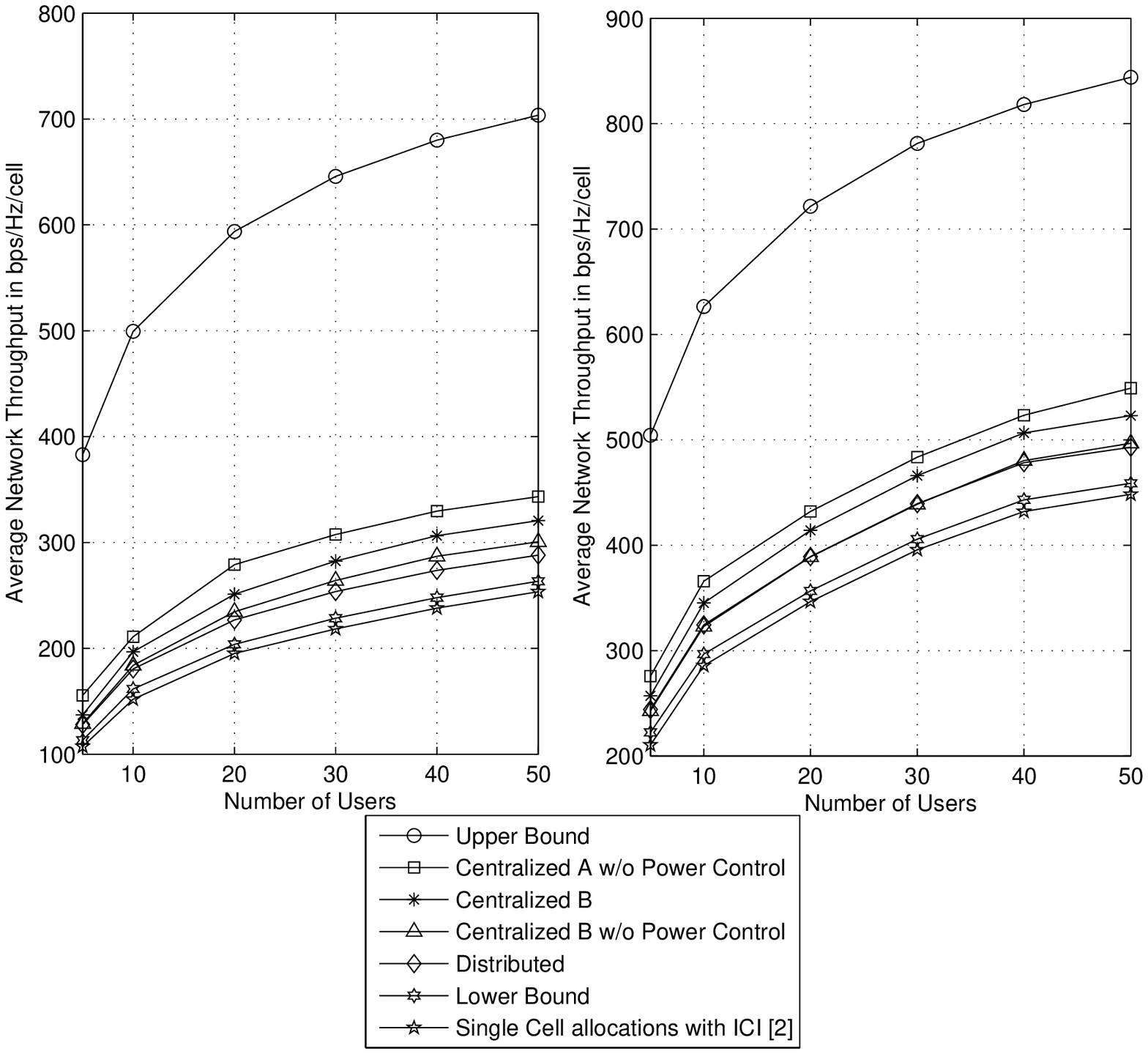}
  \caption{Comparison of all proposed schemes for $L$=7 cells, (a) Scenario A (b) Scenario B: Users are placed at 0.9 km from~ BS}
\label{Cap}
\end{figure}
\begin{figure}[H]
 \centering
\includegraphics[width=5in]{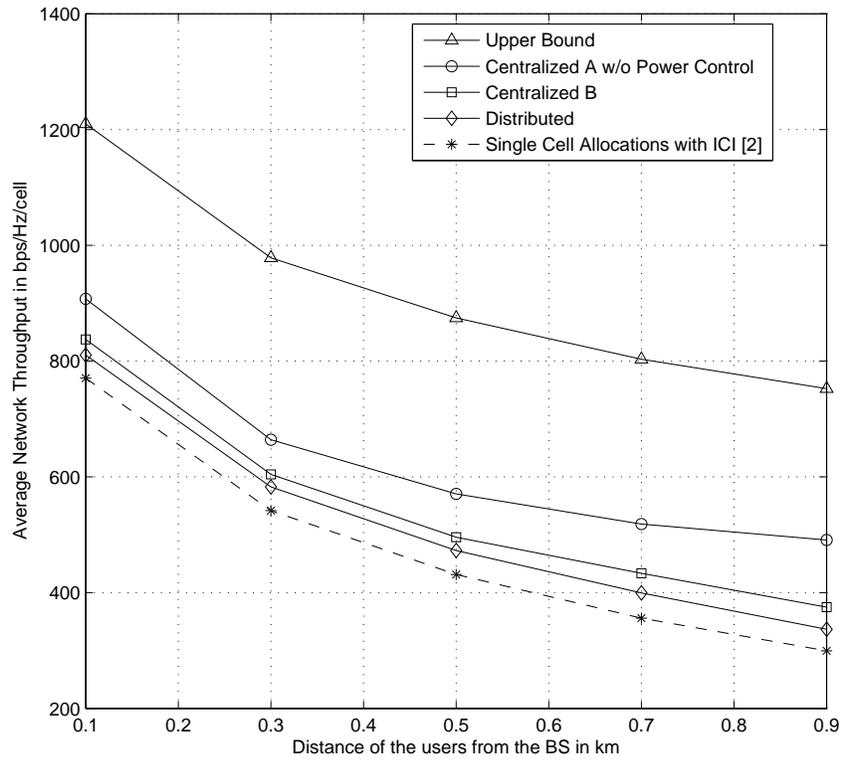}
\caption{Comparison of all proposed schemes for $L$=7 cells, $K$=100 users from cell center to cell edge}
\end{figure}
\end{document}